\documentstyle{article}
\newcommand{\nc}{\newcommand}
\nc{\be}{\begin{equation}}
\nc{\ee}{\end{equation}}
\nc{\dfrac}{\displaystyle \frac }
\author{ A. Yildiz\\
Feza G\"ursey Institute, P.O. Box 6, 81220,
Cengelk\"oy, Istanbul, Turkey}
\title{A new algebra which transmutes to the braided algebra  }

\begin{document}

\maketitle
\begin{abstract}
We find a new braided Hopf structure for the algebra
satisfied by the entries
of the braided matrix $BSL_q(2)$. A new nonbraided algebra whose coalgebra
structure is the same as the braided one is found to be a two parameter
deformed algebra. It is found that this algebra is not a comodule algebra
under adjoint coaction. However, it is shown that for a certain value of one
of the deformation parameters the braided algebra becomes a comodule algebra
under the coaction of this nonbraided algebra. The transmutation of the
nonbraided algebra to the braided one is constructed explicitly.
\end{abstract}

\vspace*{2cm}
\baselineskip=20pt
\section{Introduction}
The  covariances of  algebraic structures are of central importance
for the physical and mathematical theories constructed by using
these algebraic
structures. The investigation of quantum group covariant structures gives rise
to the introduction of braided groups
(a self contained review can be found in \cite{marev})
which are obtained from quantum
groups by a covariantization process called transmutation . Quantum groups
of function algebra type
  are not invariant under their own transformation, i.e.,
 the quantum algebra $A(R)$ is not a comodule algebra under
the adjoint coaction. However, it is possible to obtain an $A(R)$-comodule
algebra called braided algebra $B(R)$ obtained by covariantizing the
algebra of
the coacted copy keeping the coalgebra
unchanged. This covariantization process   deforms the
notion of tensor product and  leads to the generalization of the
usual Hopf algebra axioms called braided Hopf algebra axioms \cite{majid1}. 
Thus the
braided Hopf
algebras  can be used to generalize the supersymmetric structures via the
generalization of super tensor product and super Hopf algebras \cite{dunne}
or to introduce oscillators with braid statistics \cite{ali}.

In this work, we investigate the covariance properties of the algebra
satisfied by the entries of the braided matrix $BSL_q(2)$. In other words,
we investigate
if there is any Hopf algebra other the quantum algebra
the coaction of which makes the braided algebra
a comodule algebra. Since  the coacting and the coacted
algebras have the same coalgebra structure in transmutation theory, we first
investigate the general braided Hopf algebra structure. We find that there
is one more braided Hopf algebra other the one given in the literature. 
We also find
that the nonbraided Hopf algebra whose  coproduct is the same as the braided
one is a two parameter deformed Hopf algebra. This algebra, like the
quantum algebra, is not a comodule algebra under the adjoint coaction.
For a certain value of one of the deformation parameters it turns out that
the braided algebra becomes a comodule algebra under the coaction
of this one parameter nonbraided algebra. We explicitly construct the
transmutation of the noncovariant (nonbraided) algebra to the covariant
(braided) one.

\section{Prelimineries and Review}

The right coaction of a Hopf algebra $A$ on $H$ is a linear map
$\beta :H\rightarrow H\otimes A$, i.e., 

\be\label{eq:coaction9}
\beta (h)=\Sigma h^{(i)}\otimes a^{(i)},\quad h^{(i)}\in H,\quad a^{(i)}\in A
\ee

\noindent satisfying

\be \label{eq:beta}
(\beta \otimes id)\circ \beta =(id\otimes \Delta )\circ \beta \quad
(id\otimes \epsilon )\circ \beta =id.
\ee

\noindent The algebra $H$ is a right $A$-comodule algebra if the map $\beta$
is an algebra homomorphism such that

\be
\beta (h\cdot g)=\beta (h)\cdot \beta (g) \quad \forall h,g\in H.
\ee

\noindent The consistency of the algebra homomorphism requires that
$\beta (1_H)=1_H\otimes 1_A $.
 The right adjoint coaction of a Hopf algebra on itself is defined
by

\be\label{eq:adaction}
\beta (h)=\Sigma h_{(2)}\otimes S(h_{(1)})\cdot h_{(3)}
\ee

\noindent where $h_{(1)},h_{(2)},h_{(3)}$ are given by
\be
\Delta ^{2}(h)= \Sigma h_{(1)}\otimes h_{(2)}\otimes h_{(3)}.
\ee

\noindent The quantum algebra ($A(R)$) of the quantum matrix $SL_q(2)$
is generated by $a,b,c,d$ and $1$
satisfying the relations

\begin{eqnarray}\label{eq:quantumgrouprelations}
a\cdot b &=& q^{-1}b\cdot a, \nonumber \\
a\cdot c &=& q^{-1}c\cdot a,\nonumber \\
b\cdot d &=& q^{-1}d\cdot b,\nonumber \\
c\cdot d &=& q^{-1}d\cdot c, \\
b\cdot c &=& c\cdot b, \nonumber \\
a\cdot d-d\cdot a &=& (q^{-1}-q)b\cdot c \nonumber \\
a\cdot d-q^{-1}b\cdot c &=& =1 \nonumber
\end{eqnarray}

\noindent the Hopf structure of which is given by  the coproducts
\begin{eqnarray}\label{eq:cop2}
\Delta (a) &=& a\otimes a+b\otimes c, \nonumber \\
\Delta (b) &=& a\otimes b+b\otimes d,  \nonumber \\
\Delta (c) &=& c\otimes a+d\otimes c, \\
\Delta (d) &=& c\otimes b+d\otimes d  \nonumber 
\end{eqnarray}
\noindent the counits
\be
\epsilon (a)=1,\ \  \epsilon (b)=0,\ \ \epsilon (c)=0,\ \ \epsilon (d)=1,
\ee
\noindent and by the antipodes
\be
S(a)=d,\ \ S(b)=-qb,\ \ S(c)=-q^{-1}c,\ \ S(d)=a.
\ee

\noindent The $*$-algebra structure with real $q$ is given by
\be
a^{*}=d,\ \ b^{*}=-qb,\ \ c^{*}=-q^{-1}c,\ \ d^{*}=a.
\ee
\noindent The adjoint coaction of $A(R)$ on itself is then calculated to give

\begin{eqnarray}\label{eq:coaction1}
\beta (a) &=& a\otimes d\cdot a+b\otimes d\cdot c+c\otimes(-qb\cdot a)+d\otimes(-qb\cdot c) \nonumber \\
\beta (b) &=& a\otimes d\cdot b+b\otimes d\cdot d +c\otimes (-qb\cdot c)+d\otimes (-qb\cdot d)\nonumber \\
\beta (c) &=& a\otimes (-q^{-1}c\cdot a)+b\otimes (-q^{-1}c\cdot c)+c\otimes a\cdot a+d\otimes a\cdot c \\
\beta (d) &=& a\otimes (-q^{-1}c\cdot b)+b\otimes(-q^{-1}c\cdot d)+c\otimes a\cdot b+d\otimes a\cdot d. \nonumber
\end{eqnarray}
                      
\noindent It can easily be seen that these mappings do not define an algebra
homomorphism. Hence $A(R)$ itself is not an $A(R)$-comodule algebra. However,
the multiplication ($\cdot $) of the coacted copy 
(first elements in the tensor product)
is replaced by a multiplication ($\underline{\cdot } $) such that

\begin{eqnarray}\label{eq:braidedmultiplication}
a\underline{\cdot }a &=& a\cdot a\nonumber \\
a\underline{\cdot }b &=& a\cdot b\nonumber \\
a\underline{\cdot }c &=& qc\cdot a\nonumber \\
a\underline{\cdot }d &=& a\cdot d+(q-q^{-1})c\cdot b\nonumber \\
b\underline{\cdot }a &=& q^2a\cdot b\nonumber \\
b\underline{\cdot }b &=& qb\cdot b\nonumber \\
b\underline{\cdot }c &=& q^{-1}b\cdot c+(1-q^{-2})(d-a)\cdot a\nonumber \\
b\underline{\cdot }d &=& qb\cdot d-(1-q^{-2})a\cdot b\\
c\underline{\cdot }a &=& q^{-1}c\cdot a\nonumber \\
c\underline{\cdot }b &=& q^{-1}c\cdot b\nonumber \\
c\underline{\cdot }c &=& qc\cdot c\nonumber \\
c\underline{\cdot }d &=& qc\cdot d\nonumber \\
d\underline{\cdot }a &=& d\cdot a\nonumber \\
d\underline{\cdot }b &=& d\cdot b\nonumber \\
d\underline{\cdot }c &=& d\cdot c-q^{-1}(1-q^{-2})c\cdot a\nonumber \\
d\underline{\cdot }d &=& d\cdot d-q^{-1}(1-q^{-2})c\cdot b.\nonumber 
\end{eqnarray}

\noindent  Then the algebra ($B(R)$) satisfied by $a,b,c,d$ and $1$ with this
new multiplication
($\underline{\cdot }$)

\begin{eqnarray}\label{eq:braidedgrouprelations}
b\underline{\cdot }a&=& q^2a\underline{\cdot }b, \nonumber \\
c\underline{\cdot }a&=& q^{-2}a\underline{\cdot }c, \nonumber \\
a\underline{\cdot }d&=& d\underline{\cdot }a,  \\
b\underline{\cdot }c&=& c\underline{\cdot }b\ +\ (1-q^{-2})a\underline{\cdot }(d-a), \nonumber \\
d\underline{\cdot }b&=& b\underline{\cdot }d\ +\ (1-q^{-2})a\underline{\cdot }b, \nonumber \\
c\underline{\cdot }d&=& d\underline{\cdot }c\ +\ (1-q^{-2})c\underline{\cdot }a, \nonumber \\
a\underline{\cdot }d-q^2c\underline{\cdot }b &=& 1. \nonumber
\end{eqnarray}

\noindent is an $A(R)$-comodule algebra under the coaction (\ref{eq:coaction1}).
This covariantization process is called transmutation. The  transformation 
under a noncommutative  algebra ($A(R)$ for instance)  deforms
the notion of tensor  product because the transformed algebras are 
no longer  independent. The deformed tensor product is called the braided 
tensor product and denoted by $ \underline{\otimes}$. 
The coalgebra of the transmuted algebra
$B(R)$ is of the same form as the original algebra $A(R)$, i.e.,
\begin{eqnarray}\label{eq:cop1}
\underline{\Delta} (a)&=&a\underline{\otimes} a+b\underline{\otimes} c, \nonumber \\
\underline{\Delta} (b)&=&a\underline{\otimes} b+b\underline{\otimes} d,  \nonumber \\
\underline{\Delta} (c)&=&c\underline{\otimes} a+d\underline{\otimes} c, \\
\underline{\Delta} (d)&=&c\underline{\otimes} b+d\underline{\otimes} d  \nonumber \\
\underline{\epsilon} (a)&=&\underline{\epsilon} (d)=1,\ \underline{\epsilon} (b)=\underline{\epsilon} (c)=0. \nonumber
\end{eqnarray}

\noindent The antipodes
of the generators of $B(R)$ are given by

\be \label{eq:ant3}
\underline{S}(a)=q^{2}d+(1-q^{2})a,\ \ \underline{S}(b)=-q^{2}b,\ \ \underline{S}(c)=-q^{2}c,\ \ \underline{S}(d)=a 
\ee

\noindent  The coproducts define an  algebra homomorphism
in the braided tensor product 
space with the braided tensor product algebra
such that
\be
(x\underline{\otimes }y)\underline{\cdot }(w\underline{\otimes} z)=x\psi(y\underline{\otimes} z)d\quad x,y,w,z \in B.
\ee
\noindent where $\psi $ is a generalization of the permutation map in
boson algebras. The braiding relations define a nontrivial statistics
between copies of algebras. The braiding relations for the braided algebra
(\ref{eq:braidedgrouprelations}) with the  coalgebra (\ref{eq:cop1})
are given by

\begin{eqnarray}\label{eq:braiding1}
\psi (a\underline{\otimes} a)&=&a\underline{\otimes} a+(1-q^2)b\underline{\otimes} c \nonumber\\
\psi (a\underline{\otimes} b)&=&b\underline{\otimes} a \nonumber\\
\psi (a\underline{\otimes} c)&=&c\underline{\otimes} a+(1-q^2)(d-a)\underline{\otimes} c \nonumber\\
\psi (a\underline{\otimes} d)&=&d\underline{\otimes} a+(1-q^{-2})b\underline{\otimes} c \nonumber\\
\psi (b\underline{\otimes} a)&=&a\underline{\otimes} b+(1-q^2)b\underline{\otimes} (d-a) \nonumber\\
\psi (b\underline{\otimes} b)&=&q^2b\underline{\otimes} b \nonumber\\
\psi (b\underline{\otimes} c)&=&q^{-2}c\underline{\otimes} b+(1+q^2)(1-q^2)^2b\underline{\otimes}
c-(1-q^{-2})(d-a)\underline{\otimes} (d-a) \nonumber\\
\psi (b\underline{\otimes} d)&=&d\underline{\otimes} b+(1-q^{-2})b\underline{\otimes} (d-a) \nonumber\\
\psi (c\underline{\otimes} a)&=&a\underline{\otimes} c \nonumber\\
\psi (c\underline{\otimes} b)&=&q^{-2}b\underline{\otimes} c \\
\psi (c\underline{\otimes} c)&=&q^2c\underline{\otimes} c \nonumber\\
\psi (c\underline{\otimes} d)&=&d\underline{\otimes} c \nonumber\\
\psi (d\underline{\otimes} a)&=&a\underline{\otimes} d+(1-q^{-2})b\underline{\otimes} c \nonumber\\
\psi (d\underline{\otimes} b)&=&b\underline{\otimes} d \nonumber\\
\psi (d\underline{\otimes} c)&=&c\underline{\otimes} d+(1-q^{-2})(d-a)\underline{\otimes} c \nonumber\\
\psi (d\underline{\otimes} d)&=&d\underline{\otimes} d-q^{-2}(1-q^{-2})b\underline{\otimes} c. \nonumber\\
\end{eqnarray}

\noindent The central element $q^{-1}a+qd$ in quantum algebra
(which is the quantum trace in the matrix algebra)
  is not only central in the braided algebra but also 
 bosonic  in the sense that
\be
\psi((q^{-1}a+qd)\underline{\otimes}x)=x\underline{\otimes}(q^{-1}a+qd), \ \ \psi(x\underline{\otimes}(q^{-1}a+qd))=(q^{-1}a+qd)\underline{\otimes}x\ \ \forall x\in B(R).
\ee

\noindent The algebras in the braided category satisfy the braided Hopf
algebra axioms \cite{majid1}

\begin{eqnarray}
\label{eq:bhopf}
m\circ (id\otimes m)&=&m\circ (m\otimes id) \nonumber \\
m\circ (id\otimes \eta )&=&m\circ (\eta \otimes id) = id \nonumber \\
(id\otimes \Delta )\circ \Delta &=&(\Delta \otimes id)\circ \Delta \nonumber \\
(\epsilon \otimes id)\circ \Delta &=&(id\otimes \epsilon )\circ \Delta = id  \nonumber \\
m\circ (id\otimes S)\circ \Delta &=&m\circ (S\otimes id)\circ \Delta  =\eta \circ \epsilon \nonumber \\ 
\psi \circ (m \otimes id) &=& (id \otimes m)\circ (\psi \otimes id)\circ (id \otimes \psi) \nonumber \\
\psi \circ (id \otimes m) &=& (m \otimes id)\circ (id \otimes \psi)\circ  (\psi \otimes id)\nonumber \\
(id \otimes \Delta)\circ \psi &=& (\psi \otimes id)\circ (id \otimes \psi )\circ (\Delta \otimes id) \nonumber \\
(\Delta \otimes id)\circ \psi &=& (id \otimes \psi )(\psi \otimes id)\circ (id \otimes \Delta) \nonumber \\
\Delta \circ m &=& (m \otimes m)(id\otimes \psi \otimes id)\circ (\Delta \otimes \Delta ) \\
S\circ m &=& m\circ \psi \circ (S\otimes S)\nonumber \\
\Delta \circ  S &=& (S\otimes S)\circ \psi \circ \Delta  \nonumber \\
\epsilon \circ  m &=& \epsilon \otimes \epsilon \nonumber \\
(\psi \otimes id)\circ (id \otimes \psi)\circ (\psi \otimes id)&=&(id \otimes \psi)\circ (\psi \otimes id)\circ (id \otimes \psi) \nonumber
\end{eqnarray}

\noindent and the $\underline{*}$ algebras in the same category with real
deformation parameter also satisfy \cite{majid2}
\begin{eqnarray}\label{eq:*1}
\Delta \circ *  &=&  \pi \circ (* \otimes *)\circ \Delta\nonumber \\
S\circ * &=& *\circ S \\
(x\otimes y)^{*} &=& y^{*}\otimes x^{*},\ \forall x,y\in B  \nonumber
\end{eqnarray}

\noindent where we omit the underlining  for the sake of clarity and 
all mappings 
are
in the
braided sense in the axioms (\ref{eq:bhopf})-(\ref{eq:*1}). 
Note that in the $\psi\rightarrow \pi$
 limit
the axioms
(\ref{eq:bhopf})
reduce to the usual
Hopf algebra axioms. It is also possible to define more general
Hopf algebras where the counit map is no longer an algebra homomorphism
\cite{durdevic}.

The  involutions  for  (\ref{eq:braidedgrouprelations})
\be \label{eq:star1}
a^{\underline{*}}=a,\ \ b^{\underline{*}}=c,\ \ c^{\underline{*}}=d,\ \ d^{\underline{*}}=d
\ee

\noindent complete the braided Hopf $\underline{*}$-algebra structure. 

\section{A New Braided Hopf Algebra Solution}

 Expressing the braided algebra in terms of a central and bosonic 
element
\be \label{eq:p}
p\equiv q^{-2}a+d
\ee

\noindent and three other  generators $a,b,c$ makes a lot simplification in the
calculations.
The algebra  (\ref{eq:braidedgrouprelations})  can equivalently  be 
expressed as
\begin{eqnarray}\label{eq:ba1}
b\underline{\cdot }a &=& q^2a\underline{\cdot }b  \nonumber \\
a\underline{\cdot }c &=& q^2c\underline{\cdot }a \\
b\underline{\cdot }c &=& c\underline{\cdot }b-(1-q^{-4})a^2+(1-q^{-2})p\underline{\cdot }a \nonumber \\
q^{-2}a\underline{\cdot }a+a\underline{\cdot }p-q^2c\underline{\cdot } b &=& 1 \nonumber
\end{eqnarray}

\noindent The $\underline{*}$-structure (\ref{eq:star1}) implies
\be
a^{\underline{*}}=a,\ \ b^{\underline{*}}=c,\ \ c^{\underline{*}}=b,\ \ p^{\underline{*}}=p.
\ee
\noindent   We write the general forms of the coproducts

\begin{eqnarray}\label{eq:gbha}
\underline{\Delta }(a)  &=& A_1a\underline{\otimes } a+A_2b\underline{\otimes } c
+A_3c\underline{\otimes } b+A_4p\underline{\otimes } a+A_5a\underline{\otimes } p
+A_61\underline{\otimes } a  \nonumber \\
 & &+A_7a\underline{\otimes } 1+A_81\underline{\otimes } 1+A_9p\underline{\otimes }
p+A_{10}1\underline{\otimes } p+A_{11}p\underline{\otimes } 1, \nonumber \\
\underline{\Delta }(b) &=& B_1a\underline{\otimes } b+B_2b\underline{\otimes } a+B_3b\underline{\otimes } p+B_4p\underline{\otimes }
b+B_51\underline{\otimes } b+B_6b\underline{\otimes } 1, \\
\underline{\Delta }(c) &=& B_1c\underline{\otimes } a+B_2a\underline{\otimes }
c+B_3p\underline{\otimes } c+B_4c\underline{\otimes } p+B_5c
\underline{\otimes } 1+B_61\underline{\otimes } c, \nonumber \\
\underline{\Delta }(p) &=& C_1a\underline{\otimes } a+C_2b\underline{\otimes } c+C_3c
\underline{\otimes } b+C_4p\underline{\otimes } a+C_5a\underline{\otimes } p \nonumber \\
 & & +C_61\underline{\otimes } a  
 +C_7a\underline{\otimes } 1+C_81\underline{\otimes } 1+C_9p\underline{\otimes }
p+C_{10}1\underline{\otimes } p+C_{11}p\underline{\otimes } 1 \nonumber  
\end{eqnarray}

\noindent the counits
\begin{eqnarray}
\underline{\epsilon }(a) &=& e_1, \nonumber \\
\underline{\epsilon }(b) &=& \underline{\epsilon }(c)=e_2, \\
\underline{\epsilon }(p) &=& e_3  \nonumber
\end{eqnarray}

\noindent the antipodes
\begin{eqnarray}
\underline{S}(a) &=& k_1a+k_2b+k_3c+k_4p+k_5, \nonumber \\
\underline{S}(b) &=& m_1a+m_2b+m_3c+m_4p+m_5, \\
\underline{S}(c) &=& m_1a+m_2c+m_3b+m_4p+m_5, \nonumber \\
\underline{S}(p) &=& n_1a+n_2b+n_3c+n_4p+n_5   \nonumber
\end{eqnarray}

\noindent  the braidings

\begin{eqnarray}\label{eq:babraidings}
\psi (a\underline{\otimes } a) &=& g_1a\underline{\otimes } a+g_2b\underline{\otimes } c
+g_3c\underline{\otimes } b+g_4p\underline{\otimes } a+g_5a\underline{\otimes } p \nonumber \\
 & &+g_61\underline{\otimes } a+g_7a\underline{\otimes } 1+g_81\underline{\otimes } 1+g_9p\underline{\otimes }
p+g_{10}1\underline{\otimes } p+g_{11}p\underline{\otimes } 1, \nonumber \\
\psi (a\underline{\otimes } b) &=& d_1b\underline{\otimes } a+d_2a\underline{\otimes } b+d_3p\underline{\otimes }
b+d_4b\underline{\otimes } p+d_51\underline{\otimes } b+d_6b\underline{\otimes } 1, \nonumber \\
\psi (a\underline{\otimes } c) &=& f_1c\underline{\otimes }
a+f_2a\underline{\otimes } c+f_3p\underline{\otimes } c+f_4c
\underline{\otimes } p+f_51\underline{\otimes } c+f_6c\underline{\otimes } 1, \nonumber \\
\psi (b\underline{\otimes } b) &=& z_1b\underline{\otimes } b, \\
\psi (c\underline{\otimes } b) &=& c_1a\underline{\otimes } a+c_2b\underline{\otimes } c
+c_3c\underline{\otimes } b+c_4p\underline{\otimes } a+c_5a\underline{\otimes } p \nonumber \\
 & & +c_61\underline{\otimes } a+c_7a\underline{\otimes } 1+c_81\underline{\otimes } 1+c_9p\underline{\otimes }
p+c_{10}1\underline{\otimes } p+c_{11}p\underline{\otimes } 1, \nonumber \\
\psi (b\underline{\otimes } c) &=& a_1a\underline{\otimes } a+a_2b\underline{\otimes } c
+a_3c\underline{\otimes } b+a_4p\underline{\otimes } a+a_5a\underline{\otimes } p \nonumber \\
 & & +a_61\underline{\otimes } a+a_7a\underline{\otimes } 1+a_81\underline{\otimes } 1 +a_9p\underline{\otimes } p+a_{10}1\underline{\otimes }
p+a_{11}p\underline{\otimes } 1 \nonumber
\end{eqnarray}

\noindent and their $\underline{*}$-involutions to find the solutions 
for the braided Hopf algebra structure.
The
symbols with a subscript are the parameters to be determined. For the bosonic
trace ,
the  braiding of the central element  $p$ is trivial.
We substitute these general forms into the braided Hopf algebra axioms
(\ref{eq:bhopf})
 and solve the equations by using the computer programming Maple V.
We find that that there are only two solutions. 

For the first
solution the coproducts 

\begin{eqnarray}\label{eq:firsttype}
\underline{\Delta }(a) &=& a\underline{\otimes } a+b\underline{\otimes } c \nonumber \\
\underline{\Delta }(b) &=& a\underline{\otimes } b-q^{-2}b\underline{\otimes } a+b\underline{\otimes } p  \\
\underline{\Delta }(c) &=& c\underline{\otimes } a-q^{-2}a\underline{\otimes }
c+p\underline{\otimes } c \nonumber \\
\underline{\Delta }(p) &=& (q^{-2}+q^{-4})a\underline{\otimes } a+q^{-2}b\underline{\otimes } c
+c\underline{\otimes } b-q^{-2}p\underline{\otimes } a-q^{-2}a\underline{\otimes } p+p\underline{\otimes } p \nonumber
\end{eqnarray}

\noindent the counits
\begin{eqnarray}
\underline{\epsilon }(a) &=& 1,  \nonumber \\
\underline{\epsilon }(b) &=& 0 \\
\underline{\epsilon }(c) &=& 0, \nonumber \\
\underline{\epsilon }(p) &=& 1+q^{-2} \nonumber
\end{eqnarray}

\noindent the antipodes

\begin{eqnarray}
\underline{S}(a) &=& q^2(p-a), \nonumber \\
\underline{S}(b) &=& -q^2b, \\
\underline{S}(c) &=& -q^2c, \nonumber \\
\underline{S}(p) &=& p \nonumber
\end{eqnarray}

\noindent together with the braidings

\begin{eqnarray}\label{eq:firstbraidings}
\psi (a\underline{\otimes } a) &=& a\underline{\otimes } a+(1-q^2)b\underline{\otimes } c, \nonumber \\
\psi (a\underline{\otimes } b) &=& b\underline{\otimes } a, \nonumber \\
\psi (a\underline{\otimes } c) &=& c\underline{\otimes }
a+(q^2-q^{-2})a\underline{\otimes } c+(1-q^2)p\underline{\otimes } c, \nonumber \\
\psi (b\underline{\otimes } b) &=& q^2b\underline{\otimes } b, \\
\psi (c\underline{\otimes } b) &=& q^{-2}b\underline{\otimes } c, \nonumber \\
\psi (b\underline{\otimes } c) &=& (-1-q^{-2}+q^{-4}+q^{-6})a\underline{\otimes }
a+(q^2-1-q^{-2}+q^{-4})b\underline{\otimes } c \nonumber \\
 & & +q^{-2}c\underline{\otimes } b+  
(1-q^{-4})a\underline{\otimes } p+(1-q^{-4})p\underline{\otimes } a+(q^{-2}-1)p\underline{\otimes } p \nonumber
\end{eqnarray}

\noindent and their $\underline{*}$-involutions define a braided Hopf
$\underline{*}$-algebra.
For the second solution the
coproducts

\begin{eqnarray}\label{eq:secondtype}
\underline{\Delta }(a) &=& a\underline{\otimes } a+q^4c\underline{\otimes } b, \nonumber \\
\underline{\Delta }(b) &=& -q^2a\underline{\otimes } b+b\underline{\otimes } a+q^2p\underline{\otimes } b, \\
\underline{\Delta }(c) &=& -q^2c\underline{\otimes } a+a\underline{\otimes }
c+q^2c\underline{\otimes } p, \nonumber \\
\underline{\Delta }(p) &=& (1+q^2)a\underline{\otimes } a+q^2b\underline{\otimes } c+q^4c\underline{\otimes } b
  -q^2p\underline{\otimes } a-q^2a\underline{\otimes }
p+q^2p\underline{\otimes } p \nonumber
\end{eqnarray}

\noindent the counits

\begin{eqnarray}
\underline{\epsilon }(a) &=& 1, \nonumber \\
\underline{\epsilon }(b) &=& 0, \\
\underline{\epsilon }(c) &=& 0, \nonumber \\
\underline{\epsilon }(p) &=& 1+q^{-2} \nonumber
\end{eqnarray}

\noindent and the antipodes

\begin{eqnarray}
\underline{S}(a) &=& -q^{-2}a+p, \nonumber \\
\underline{S}(b) &=& -q^{-2}b, \\
\underline{S}(c) &=& -q^{-2}c, \nonumber \\
\underline{S}(p) &=& p \nonumber
\end{eqnarray}

\noindent together with the braidings

\begin{eqnarray}
\psi (a\underline{\otimes } a) &=& a\underline{\otimes } a+(q^4-q^2)c\underline{\otimes } b, \nonumber \\
\psi (a\underline{\otimes } b) &=& b\underline{\otimes } a+(q^{-2}-q^2)a\underline{\otimes }
b+(q^2-1)p\underline{\otimes } b, \nonumber \\
\psi (a\underline{\otimes } c) &=& c\underline{\otimes } a \nonumber \\
\psi (b\underline{\otimes } b) &=& q^{-2}b\underline{\otimes } b, \\
\psi (b\underline{\otimes } c) &=& q^2c\underline{\otimes } b, \nonumber \\
\psi (c\underline{\otimes } b) &=& (q^2+1-q^{-2}-q^{-4})a\underline{\otimes }
a+q^2b\underline{\otimes } c+(q^4-q^2+q^{-2}-1)c\underline{\otimes } b \nonumber \\
 & & +(q^{-2}-q^2)a\underline{\otimes } p+(q^{-2}-q^2)p\underline{\otimes } a+(q^2-1)p\underline{\otimes } p \nonumber
\end{eqnarray}

\noindent and their $\underline{*}$-involutions makes $B(R)$
a braided Hopf $\underline{*}$-algebra.
Note that
in the $q\rightarrow 1$ limit, not only the algebra becomes commutative but
also the braided tensor product becomes the ordinary tensor product, and
we obtain the bosonic statistics for both braided Hopf algebra solutions.
When we express the solutions in terms of the original generators of the 
algebra, i.e., in terms of $a,b,c,d$ by using (\ref{eq:p}) we see that
that the first solution is completely equivalent to the braided Hopf algebra
given in the literature which we give in (\ref{eq:cop1})-(\ref{eq:braiding1}).
But the second solution with  
the coproducts
\begin{eqnarray}\label{eq:cop3}
\underline{\Delta }(a) &=& a\underline{\otimes } a+q^4c\underline{\otimes } b, \nonumber \\ 
\underline{\Delta }(b) &=& (1-q^2)a\underline{\otimes } b+b\underline{\otimes } a+q^2d\underline{\otimes } b, \nonumber \\ 
\underline{\Delta }(c) &=& (1-q^2)c\underline{\otimes } a+a\underline{\otimes } c+q^2c\underline{\otimes } d,  \\
\underline{\Delta }(d) &=& (q^2-1)a\underline{\otimes } a+(q^4-q^2)c\underline{\otimes } b+q^2b\underline{\otimes } c 
            +(1-q^2)a\underline{\otimes } d+(1-q^2)d\underline{\otimes } a+q^2d\underline{\otimes } d \nonumber  
\end{eqnarray}
\noindent  counits
\be \label{eq:counit3}
\underline{\epsilon }(a)=\underline{\epsilon }(d)=1,\ \ \  \underline{\epsilon }(b)= \underline{\epsilon }(c)=0, 
\ee
\noindent antipodes
\be
\underline{S}(a)=d,\ \ \underline{S}(b)=-q^{-2}b,\ \ \underline{S}(c)=-q^{-2},\ \ \underline{S}(d)=q^2d+(1-q^2)a
\ee
\noindent and braidings
\begin{eqnarray}\label{eq:bra2}
\psi (a\underline{\otimes } a)&=&a\underline{\otimes } a+(q^4-q^2)c\underline{\otimes } b, \nonumber \\ 
\psi (a\underline{\otimes } b)&=&b\underline{\otimes } a+ (1-q^2)a\underline{\otimes } b+(q^2-1)d\underline{\otimes } b,  \nonumber \\ 
\psi (a\underline{\otimes } c)&=&c\underline{\otimes } a ,  \nonumber \\ 
\psi (a\underline{\otimes } d)&=&d\underline{\otimes } a+(1-q^2)c\underline{\otimes } b ,  \nonumber \\ 
\psi (b\underline{\otimes } a)&=&a\underline{\otimes } b , \nonumber \\ 
\psi (b\underline{\otimes } b)&=&q^{-2}b\underline{\otimes } b, \nonumber \\ 
\psi (b\underline{\otimes } c)&=&q^2c\underline{\otimes } b,  \nonumber \\ 
\psi (b\underline{\otimes } d)&=&d\underline{\otimes } b,  \\
\psi (c\underline{\otimes } a)&=&a\underline{\otimes } c +(1-q^2)c\underline{\otimes } a+(q^2-1)c\underline{\otimes } d, \nonumber \\ 
\psi (c\underline{\otimes } b)&=&(q^2-1)a\underline{\otimes } a+q^2b\underline{\otimes } c+(q^4-q^2+q^{-2}-1)c\underline{\otimes } b \nonumber \\ 
                 & &+(1-q^2)a\underline{\otimes } d+(1-q^2)d\underline{\otimes } a+(q^2-1)d\underline{\otimes } d,   \nonumber \\ 
\psi (c\underline{\otimes } c)&=&q^{-2}c\underline{\otimes } c, \nonumber \\ 
\psi (c\underline{\otimes } d)&=&d\underline{\otimes } c+(1-q^{-2})c\underline{\otimes } a+(q^{-2}-1)c\underline{\otimes }  d, \nonumber \\ 
\psi (d\underline{\otimes } a)&=&a\underline{\otimes } d+(1-q^2)c\underline{\otimes } b, \nonumber \\ 
\psi (d\underline{\otimes } b)&=&b\underline{\otimes } d+(1-q^{-2})a\underline{\otimes } b+(q^{-2}-1)d\underline{\otimes } b, \nonumber \\ 
\psi (d\underline{\otimes } c)&=&c\underline{\otimes } d, \nonumber \\ 
\psi (d\underline{\otimes } d)&=&d\underline{\otimes } d+(1-q^{-2})c\underline{\otimes } b \nonumber  
\end{eqnarray}

\noindent defines a different braided Hopf algebra.

\section{A New Algebra and its Transmutation}

We find that the nonbraided Hopf algebra generated by the generators
$a,b,c,d$ and $1$
whose coalgebra is of the form (\ref{eq:cop3})-(\ref{eq:counit3}) is a two
parameter ($q,r$) deformed Hopf algebra. The algebra part is found to satisfy 

\begin{eqnarray}
a\cdot b &=& rb\cdot a, \nonumber \\ 
a\cdot c &=& rc\cdot a, \nonumber \\ 
b\cdot c &=& c\cdot b,  \\
b\cdot d &=& rd\cdot b+(q^{-2}-1)(r^2-1)b\cdot a \nonumber \\ 
c\cdot d &=& rd\cdot c+(q^{-2}-1)(r^2-1)c\cdot a, \nonumber \\ 
a\cdot d-d\cdot a &=& (r-r^{-1})q^2b\cdot c, \nonumber \\ 
q^2d\cdot a+(1-q^2)a\cdot a-r^{-1}q^4c\cdot b &=& 1.\nonumber 
\end{eqnarray}

\noindent Finally the antipodes 
\be
S(a)=(1-q^2)a+q^2d,\ S(b)=-rb,\ S(c)=-r^{-1}c,\ S(d)=(2-q^2)a+(q^2-1)d
\ee
\noindent completes Hopf algebra. For the $*$-algebra structure we find
\be
a^{*}=(1-q^2)a+q^2d,\ \ b^{*}=-r^{-1}c,\ \ c^{*}=-rb,\ \
d^{*}=(2-q^2)a+(q^2-1)d.
\ee
The adjoint coaction of this Hopf algebra on itself calculated by using
(\ref{eq:adaction})

\begin{eqnarray}\label{eq:coaction}
\beta (a) &=& a\otimes [(1-q^2)a\cdot a+q^2d\cdot a-r^{-1}q^4(1-q^2)c\cdot b]+b\otimes [-r^{-1}q^4c\cdot a] \nonumber \\
          & & +c\otimes [q^4(1-q^2)a\cdot b+q^6d\cdot b]+d\otimes [-r^{-1}q^6c\cdot b],\nonumber \\
\beta (b) &=& a\otimes [-q^2a\cdot b]+b\otimes [a^2]+c\otimes [-q^4rb\cdot b]+d\otimes [q^2a\cdot b] \\
\beta (c) &=& a\otimes [q^4d\cdot c+q^2(1-q^2)a\cdot c]+c\otimes [(1-q^2)^2a\cdot a+q^2(1-q^2)a\cdot d+q^2(1-q^2)d\cdot a+q^4d\cdot d]\nonumber \\
          & & b\otimes [-q^4r^{-1}c\cdot c]+d\otimes [q^2(q^2-1)r^{-1}c\cdot a+q^4r^{-1}c\cdot d]\nonumber \\
\beta (d) &=& a\otimes [(q^6-q^4)r^{-1}c\cdot b-q^4rc\cdot b]+c\otimes [(q^6-q^4-q^4r^2)d\cdot b+(q^4-q^2)(r^2-q^2+1)a\cdot b]\nonumber \\
          & & b\otimes [(q^2(1-q^2)r^{-1}+q^2r)c\cdot a]+d\otimes [(1-q^2)a\cdot a+q^2a\cdot d-(q^6-q^4)r^{-1}c\cdot b]\nonumber
\end{eqnarray}

\noindent does not define an algebra homomorphism. However when the
coacted copy satisfies  the braided algebra (\ref{eq:braidedgrouprelations}),
the transformations (\ref{eq:coaction}) define an algebra homomorphism for 
$r=q$.
It can be shown that the braiding of the transformation in the braided 
tensor product space is equal to the transformation of the braiding on the 
same space, i.e., 
\be
\psi (\beta (x\underline{\otimes }y))=\beta (\psi (x\underline{\otimes }y))\ \ \ \forall x,y\in B(R)
\ee
\noindent where the braidings are given by (\ref{eq:bra2}). It can also be 
shown that the braided  algebra is a comodule algebra under the 
adjoint coaction for both of the braided Hopf algebra solutions.

The transmuted
multiplication ($\underline{\cdot }$)
in terms of the multiplication ($\cdot$) of the coacting nonbraided algebra 

\begin{eqnarray}
a\underline{\cdot } a &=&a\cdot a \nonumber \\                      
a\underline{\cdot } b &=&q^{-1}b\cdot a \nonumber \\ 
a\underline{\cdot } c &=&a\cdot c  \nonumber \\ 
a\underline{\cdot } d &=&a\cdot d+(q-q^3)b\cdot c  \nonumber \\ 
b\underline{\cdot } a &=&qb\cdot a   \nonumber \\ 
b\underline{\cdot } b &=&q^{-1}b\cdot b \nonumber \\ 
b\underline{\cdot } c &=&qb\cdot c  \\
b\underline{\cdot } d &=&q^{-1}b\cdot d+q(1-q^{-2})^2b\cdot a \nonumber \\ 
c\underline{\cdot } a &=&q^{-2}a\cdot c \nonumber \\ 
c\underline{\cdot } b &=&qc\cdot b+(q^{-2}-1)(d-a)\cdot a \nonumber \\ 
c\underline{\cdot } c &=&q^{-1}c\cdot c \nonumber \\ 
c\underline{\cdot } d &=&q^{-1}c\cdot d+(1-q^{-2})a\cdot c \nonumber \\ 
d\underline{\cdot } a &=&d\cdot a \nonumber \\ 
d\underline{\cdot } b &=&d\cdot b+(q^{-2}-q^{-4})a\cdot b \nonumber \\ 
d\underline{\cdot } c &=&d\cdot c \nonumber \\ 
d\underline{\cdot } d &=&d\cdot d+(q-q^{-1})b\cdot c \nonumber  
\end{eqnarray}

\noindent completes the transmutation process.

\end{document}